\documentclass[12pt]{amsart}

\usepackage{amsfonts,newlfont,latexsym}
\usepackage{euscript}
\usepackage{amssymb,amscd}

\usepackage{amsmath}

      \theoremstyle{plain}
      \newtheorem{theorem}{Theorem}[section]
      \newtheorem{lemma}[theorem]{Lemma}
      \newtheorem{corollary}[theorem]{Corollary}
      \newtheorem{proposition}[theorem]{Proposition} 
    
      \newtheorem{conjecture}[theorem]{Conjecture}
      
      \theoremstyle{definition}
        \newtheorem{remark}{Remark}

      \makeatletter
      \def\@setcopyright{}
      \def\serieslogo@{}
      \makeatother
\newcommand{\N}{\mathbb N}

\newcommand{\R}{\mathbb R}
\newcommand{\Rk}{\mathbb R^k}

\newcommand{\Z}{\mathbb Z}
\newcommand{\Zk}{\mathbb Z^k}

\newcommand{\T}{\mathbb T}
\newcommand{\Tk}{\mathbb T^k}

\newcommand{\tk}{\mathbb T ^{k+1}}

\newcommand{\Ce}{C^{1+\epsilon}}

\newcommand{\e}{\varepsilon}
\newcommand{\la}{\lambda}
\newcommand{\La}{\Lambda}

\def \e{\varepsilon}

\def\dist{\text{dist}}

\def \a{\alpha}
\def \am{{\alpha (\mathbf m)}}
\def \an{{\alpha (\mathbf n)}}
\def \at{{\alpha (\mathbf t)}}
\def \as{{\alpha (\mathbf s)}}
\def \ao{{\alpha _0}}
\def \aom{{\alpha _0 (\mathbf m)}}
\def \aon{{\alpha _0 (\mathbf n)}}
\def \aot{{\alpha _0 (\mathbf t)}}

\def \m{\mathbf m}
\def \n{\mathbf n}
\def \t{\mathbf t}
\def \s{\mathbf s}

\def \h{\mathbf h}

\def \wl{\tilde W}
\def \w{\tilde {\cal W}}
\def \wol{W}
\def \wo{\cal W}
\def \bwr{\tilde B_r}

\def \bwo{B}
\def \mw{\mu ^{\tilde {\cal W}}}

\def \de{D^E}
\def \df{Jf}
\def \dfn{Jf^n}      

\def\proof{{\bf Proof.}\; }
\def\QED{\hfill\hfill{\square}}

\numberwithin{equation}{section}
\begin{document}

\author[Boris Kalinin and Anatole Katok]{Boris Kalinin $^1)$ and Anatole Katok $^2)$ }
\address {University of South Alabama, Mobile, AL}
\email{kalinin@jaguar1.usouthal.edu}

  \address{The Pennsylvania State University, University Park, PA}

\email{katok\_a@math.psu.edu}

\title[Rigidity beyond uniform hyperbolicity]{Measure rigidity beyond uniform hyperbolicity: Invariant Measures for Cartan actions on Tori}

%%%%%%%%%%%%%%   Abstract   %%%%%%%%%%%%%%%%

\begin{abstract}
We prove that every smooth action $\a$ of  $\Zk,\,\,k\ge 2$, on the $(k+1)$-dimensional 
torus homotopic to an action $\ao$  by  hyperbolic linear maps preserves 
an absolutely continuous measure. This is a first known result concerning abelian groups 
of diffeomorphisms where existence of an invariant geometric  structure is  obtained from homotopy data. 

We also show that both ergodic and geometric properties of such a measure are very close to 
the corresponding properties of the Lebesgue measure with respect to the linear action $\ao$.

\end{abstract}

%%%%%%%%%%%%%%  End Abstract   %%%%%%%%%%%%%%%

\date{\today}

\thanks{$^1)$ This work has been partially supported by NSF grant DMS-0140513 and by the Center for Dynamics and Geometry at Penn State.}

\thanks{$^2)$  Partially supported by NSF grant DMS-0505539 }

\maketitle

%%%%%%%%%%%%%%%%% Begin Text %%%%%%%%%%%%%%%%%%

\section{Introduction}
\subsection{Measure rigidity and hyperbolicity}
It is well-known that in classical dynamical systems, i.e. smooth actions of $\Z$ or $\R$,
non-trivial recurrence combined with  some kind of hyperbolic behavior  produces a rich 
variety of invariant measures (see for example, \cite{KH} and \cite{KM} for the uniformly 
and non-uniformly hyperbolic situations correspondingly). On the other hand, invariant 
measures for actions of {\em higher rank} abelian groups tend to be scarce. This was 
first  noticed by Furstenberg \cite{F} who posed the still open problem of describing all 
ergodic measures on the circle invariant with respect to  multiplications by  2 and  by 3. 
Great progress  has been made in characterizing invariant measures with positive entropy 
for algebraic  actions of higher rank abelian groups; for the   measure rigidity  results for 
actions by automorphisms or endomorphisms of  a torus see 
\cite{R, KS3, KS4, KaK1, KaK2, KaSp, EL}.

For the background on algebraic, arithmetic, and ergodic properties of $\Zk$ actions 
by automorphisms of the torus we refer to \cite{KKS}. Recall  that an action of $\Zk$ 
on $\tk$, $k \ge2$, by  automorphisms which are  ergodic with respect to the Lebesgue measure 
is called a (linear) {\em Cartan action}. Every element of a Cartan action other than 
identity is hyperbolic and  has distinct real eigenvalues, and the centralizer of a Cartan 
action in the groups of automorphisms of the torus is a finite extension of the action 
itself (\cite[Section 4.1]{KKS}).

A geometric approach to measure rigidity  was introduced in \cite{KS3}. It is  based 
on the study of conditional measures on various invariant foliations for the action. 
Broadly speaking, there are three essential  tools or methods  within this approach which we list 
in order of their chronological  appearance:
\smallskip

\begin{enumerate}
\item\label{Lyapunovplanes} Geometry of  Lyapunov exponents and  derivative objects, 
in particular Weyl chambers \cite{KS3, KS4, KaK1, KaSp}. 
\footnote{$\,$ See in particular \cite[Section 2.2]{KaK1} for a  
down-to-earth proof of rigidity of positive entropy invariant measures for linear  Cartan actions. A reader 
unfamiliar with  measure rigidity may  look  al that section first to get an idea of 
the  basic arguments  we   generalize in the present paper.}
\item The non-commutativity and specific  commutation relations between  various 
invariant foliations \cite{EK1, EK2, EKL}. 
\footnote{$\,$ This   method was first outlined at the end of \cite{KS3}; notice that it is not relevant 
for the actions on the torus since in this case all foliations commute.}
\item Diophantine properties of global recurrence \cite{EL}.
\end{enumerate}
\smallskip

In this paper  we make the first step   in  extending  measure rigidity from algebraic actions 
to  the general non-uniformly  hyperbolic case, i.e. to positive entropy ergodic invariant 
measures  for actions of higher rank abelian groups all of whose Lyapunov characteristic exponents 
do not vanish. Such measures are usually called {\em hyperbolic measures}. 
The theory of hyperbolic
measures for smooth actions of higher rank abelian groups is described in Part II of \cite{KaK1}. 
In Sections~\ref{section-prelim} and \ref{sectionPesin} we  will  briefly mention certain key elements 
of that theory relevant for the  specific  situation considered in this paper.   
 
 In this paper  we use a counterpart of  the method (\ref{Lyapunovplanes}) above.  
 We will discuss the scope of this method, difficulties which  appear for its extensions, and applications of properly modified versions of other methods to various non-uniformly hyperbolic situations in a subsequent paper.
 \medskip
 
\noindent{\bf Acknowledgement.}  We would like to thank Omri Sarig who carefully read the paper and made a number of valuable comments which helped to clarify several points  in the  proofs and improve presentation.
 
\subsection{Formulation of results}

\begin{theorem} \label{main}
Any action $\a$  of $\Zk, \,\,k\ge 2$, by $C^{1+\epsilon},\,\epsilon>0$, diffeomorphisms of 
$\tk$, which is homotopic to a  linear Cartan action $\a _0$,  has an ergodic absolutely 
continuous invariant measure.
\end{theorem}
\vskip0.3cm

The connection between  invariant measures of $\a$ and those of $\ao$ is established using
the following well-known result  whose proof we include for the sake of completeness.

\begin{lemma}\label{lemmaFranks}
There is a  unique surjective continuous map $h : \tk \to \tk$ homotopic to identity  such that
$$
h \circ \a=\a_0 \circ h.
$$
\end{lemma}
\proof
 Consider an element $\m\in \Zk\setminus\{ 0\}$.
By a theorem of Franks (\cite[Theorem 2.6.1]{KH}), there exists a unique continuous map 
$h : \tk \to \tk$ that is homotopic to identity and satisfies
\begin{equation}
h \circ \am=\aom \circ h. \label{1}
\end{equation}
For any other element $\m' \in \Zk$ consider the map 
 \begin{equation}\label{2}
  h'=\ao (-\m') \circ h \circ \a (\m')  
\end{equation}
Using commutativity of both actions $\a$ and $\ao$ as well as \eqref{1} we obtain
 $$
  h' \circ \am=\ao (-\m') \circ h \circ \a (\m') \circ \am =\ao (-\m') \circ h \circ \am \circ \a(\m')= 
 $$
 $$
 \ao (-\m') \circ \aom \circ h \circ \a(\m')=\aom \circ \ao (-\m')\circ h \circ \a(\m')=\aom \circ h',
 $$
i.e. $h'$ satisfies \eqref{1}. Since it is also homotopic
to identity, the uniqueness of $h$ forces $h=h'$. Then \eqref{2} implies
$$
h \circ \a (\m') =\ao (\m') \circ h
$$
and thus $h$ intertwines the actions $\a$ and $\ao$. $\QED$

Another way of stating Lemma~\ref{lemmaFranks} is that  the algebraic 
action $\a_0$ is a topological factor of the action $\a$ or, equivalently, 
$\a$ is an extension of $\a_0$.

\begin{remark} If  the action $\a$ is Anosov, i.e. if $\am$ is an Anosov diffeomorphism 
for some $\m$, then the map $h$ is invertible, and  both $h$ and $h^{-1}$ are H\"older
\cite[Theorems 18.6.1 and 19.1.2]{KH}. This implies various rigidity results  for $\Zk$ 
Anosov actions on the torus.  
\begin{itemize}
\item
For example, if $\ao$ is a linear $\Zk$ action  on a torus 
which contains a $\Z^2$ subaction all of whose elements other than  identity are ergodic,
then any Anosov action $\a$ homotopic to $\ao$ preserves a smooth measure.
This follows from rigidity of H\"older cocycles over $\ao$ and hence over $\a$ applied to the logarithm of the Jacobian for $\a$.
\item  
For those cases when  positive entropy ergodic invariant measure for $\ao$  is unique 
\cite{EL}  the same is true for $\a$.
\end{itemize}
\end{remark}

Consider  the set of all Borel probability measures $\nu$ on $\tk$ such that 
$(h)_*\nu=\lambda$, where $\lambda$ is  Lebesgue measure on $\tk$. 
This set is  convex, weak* compact, and $\a$ invariant. Hence by Tychonoff theorem 
it contains a nonempty subset $\mathcal M$ of measures  invariant under $\a$. 
Since $\a_0$ is ergodic with respect to $\lambda$, almost every ergodic component 
of an $\a$-invariant measure $\nu\in\mathcal M$ also belongs to $\mathcal M$. 
Let $\mu$ be such an ergodic measure.

Theorem~\ref{main}  follows immediately from Lemma~\ref{lemmaFranks} 
and the following theorem, which is the first principal technical result of the present paper. 

\begin{theorem}\label{thm-projection} Any ergodic $\a$-invariant measure $\mu$ 
such that $(h)_*\mu=\lambda$, where $h$ is the semiconjugacy from 
Lemma~\ref{lemmaFranks}, is  absolutely continuous.
\end{theorem}
 Since any $\a$-invariant measure whose ergodic components are absolutely continuous is itself  absolutely continuous  we obtain the following.

\begin{corollary}
Every measure $\nu\in\mathcal M$ is absolutely continuous and has no more 
than countably many ergodic components. Hence $\mathcal M$ contains at most countably many ergodic measures.
\end {corollary}

We believe that a much stronger statement  should be true.
%\foot{true for single map?}

\begin{conjecture} 
$\,$

1. The set $\mathcal M$ consists of a single measure.   

2. The semiconjugacy $h$ is a measurable isomorphism  between actions  $\a$ and $\ao$. 

\end{conjecture}

While part 1 of this conjecture, or even  finiteness of ergodic measures in $\mathcal  M$, 
remains open, we prove an only slightly weaker version of part 2 for ergodic measures.

\begin{theorem}\label{thm-finitecover} For any ergodic measure $\mu\in\mathcal M$ the 
semiconjugacy $h$ is finite-to-one in the following sense. There is an $\a$-invariant set 
$A$ of full measure $\mu$  such that for $\lambda$ almost every $x\in\tk$, $A\cap h^{-1}(\{x\})$
consists of equal number $s$ of points and the conditional measure induced by $\mu$ assigns 
every point in  $A\cap h^{-1}(\{x\})$ equal measure $1/s$.
\end{theorem}

Recall that  Lyapunov characteristic exponents of the linear action $\ao$ are independent of 
an invariant measure and are equal to the  logarithms of  the absolute values of the eigenvalues.  They all have multiplicity one and no two of them are proportional. 

\begin{theorem}\label{thm-exponents} Lyapunov characteristic exponents of the action 
$\a$ with respect to any ergodic measure $\mu\in\mathcal M$ are equal to the Lyapunov 
characteristic exponents of the  action $\ao$. 
\end{theorem}

Either of the last two theorems immediately implies the following. 

\begin{corollary}
The entropy function of $\a$ with respect to any measure $\nu\in\mathcal M$ is 
the same as the entropy function of $\ao$ with respect to  Lebesgue measure,
i.e. for any measure $\nu\in\mathcal M$ and any $\m \in \Zk$
$$
\h_{\nu}(\a(\m)) = \h_{\lambda}(\a_0(\m)).
$$
\end {corollary}

%\foot{formulate corollaries about factors, centralizers and joinings which follow from corresponding results for Cartan actions}
Since every element of $\ao$ other than identity is Bernoulli with respect to the
Lebesgue measure, Theorem~\ref{thm-finitecover} 
also implies that every element  of $\a$ is Bernoulli up to a finite permutation. 

\begin{corollary}\label{cor-Bernoulli} There exist a partition  of a set $A$  of full 
measure $\mu$ into finitely 
many sets $A_1,\dots , A_m$ of equal measure such that every element of $\a$ permutes  
these sets. Furthermore, there is a subgroup of finite index $\Gamma\subset \Zk$ such that 
for any $\gamma\in\Gamma$ other than identity $\a(\gamma)A_i=A_i,\,\,i=1,\dots, m$, 
and the restriction of $\a(\gamma)$ to each set $A_i$ is Bernoulli. 

In particular, if all non-identity elements of $\a$ are ergodic then they are Bernoulli. 
\end{corollary}

\begin{remark} Since  by Lemma~\ref{Weyl chambers} the  measure $\mu$  is hyperbolic   Corollary~\ref{cor-Bernoulli}  follows directly from Theorem~\ref{thm-projection} and the  classical result of Pesin \cite{P} which states that  any  ergodic hyperbolic absolutely continuous  measure for a diffeomorphism  is Bernoulli up to a finite permutation.
\end{remark}

\begin{remark}
Theorem~\ref{main} for $k=2$ was announced in \cite{KaK1} as Theorem 8.2.  
The proof  in the present paper follows a path different from the one outlined in \cite{KaK1}. 
At the moment we  do not have a complete proof of Theorem~\ref{main}  which follows the 
scheme outlined in \cite{KaK1}.   Notice that  our Theorems~\ref{thm-finitecover} and 
\ref{thm-exponents}  and their corollaries  give considerably more detailed information 
about the structure of absolutely continuous invariant measures for actions homotopic 
a to linear Cartan action than what  follows from the results announced in  \cite{KaK1}.
\end{remark}

\section{Lyapunov exponents, Weyl chambers, and invariant ``foliations'' for $\a$}

\subsection{Preliminaries}\label{section-prelim}

\subsubsection{Entropy}
Since $h_*\mu=\lambda$ the measure-theoretic entropy $\h_{\mu}$ satisfies
$$
\h_{\mu}(\a(\m))\ge \h_{\lambda}(\a_0(\m))\ge\max_{1\le i\le k+1}|\log|\rho_i(\m)|\,\,|,
$$
where $\rho_i(\m),\,\,i=1,\dots,k+1$ are  the
eigenvalues of the matrix $\aom$.

Since every element of $\ao$ other than identity is hyperbolic this implies, in particular, that
\smallskip

$(\mathcal E)$ {\em The entropies $\h_{\mu}(\a(\m))$ for
all $\m\in \Zk\setminus\{ 0\}$ are uniformly bounded away from zero.} 

\subsubsection{Lyapunov exponents}
The linear functionals on $\Zk$,  $\chi_i=\log |\rho_i|$, $i=1,\dots,k+1$  are the {\em Lyapunov 
characteristic exponents} of the linear action $\ao$ which are independent of an 
invariant measure.  See \cite[Section 1.2]{KaK1} for  the definitions and  discussion 
of Lyapunov characteristic exponents, related notions (Lyapunov hyperplanes, Weyl 
chambers, {\em etc}.) and suspensions in this setting.  We will use this  material without 
further references.

The following property of linear Cartan actions will  play an important role in our considerations, in particular  in Section~\ref{ergodicitysingular}
\smallskip

$(\mathcal C)$ {\em For every $i\in\{1,\dots,k+1\}$ there exists an element $\m\in\Zk$ 
such that $\chi_i(\m)<0$ and $\chi_j(\m)>0$ for all $j\neq i$}. (The same inequalities 
hold for any other element $\m '$ in the Weyl chamber of $\m$.) 
\medskip

Corresponding notions in a general setting, which includes that of $\Zk$ actions by 
measure preserving diffeomorphisms of smooth manifolds, are defined and discussed 
in Sections 5.1 and 5.2 of the same paper \cite{KaK1}. We will also use those notions without 
special references.  

Let $\tilde\chi_i,\,\,i=1,\dots , k+1$, be the Lyapunov characteristic exponents of the action 
$\a$, listed with their multiplicities if necessary. We will eventually show that in our setting
the exponents can be properly numbered so that  $\tilde\chi_i=\chi_i,\,\,i=1,\dots ,k+1$
(see Section~\ref{exponents-proof}). 

As the first step in this  direction  we will show in Section~\ref{weyl} that exponents  
for $\a$ can be numbered in such a way that they become  {\em proportional} 
to $\chi_i$  with  positive  scalar coefficients. 

\subsubsection{Suspensions}
Although the Lyapunov characteristic exponents for a $\Zk$ action are defined as linear 
functionals on $\Zk$, it seems natural to extend them  to $\Rk$. For example, Lyapunov 
hyperplanes (the kernels of the functionals) may be irrational and hence ``invisible'' within 
$\Zk$. It is natural to try to  construct an $\Rk$ action for which the extensions of the 
exponents from $\Zk$ will provide the non-trivial exponents. 

This is given by the suspension construction  which associates  to  a given $\Zk$ action 
on a space $N$ an $\Rk$ action on a bundle over $\Tk$ with fiber $N$. 
The topological 
type of the suspension space depends only on the homotopy type of the $\Zk$ action. 
In particular, the suspension spaces for $\ao$ and $\a$ are homeomorphic. There is a natural
correspondence between the invariant measures, Lyapunov exponents, Lyapunov 
distributions, stable and unstable manifolds, etc. for the original $\Zk$ action and its
suspension. Naturally, the suspension has additional $k$ Lyapunov exponents 
corresponding to the orbit directions which are identically equal to zero. In our setting, 
the semiconjugacy between $\a$ and $\ao$ naturally extends to the suspension. 
The extended semiconjugacy is smooth along the suspension orbits and reduces to 
the original semiconjugacy in the fiber over the origin in $\T^{k}$.
 
At various stages of the subsequent arguments it will be more convenient to deal either with 
the original  actions $\a$ and $\ao$ on $\tk$ or with their suspensions.
So we will take a certain liberty with the notations and {\em will use the same notations
for the corresponding objects}, i.e. $\a$ and $\ao$ for the suspension actions,  $\tilde\chi_i$ 
and $\chi_i$ for the Lyapunov exponents etc, modifying the notations when necessary, as 
in $\am$ for $\m\in\Zk$ and $\at$ for $\t\in\Rk$.  

\subsection{Pesin sets and invariant manifolds}\label{sectionPesin}

We will use the standard material on invariant  manifolds corresponding to the negative 
and positive Lyapunov exponents (stable and unstable manifolds)  for $C^{1+\epsilon}$
measure preserving diffeomorphisms of compact manifolds. See for example 
\cite[Chapter 4]{BP}. It is customary to use words ``distributions'' and ``foliations'' in this setting although in fact we are dealing with measurable  families of tangent spaces defined almost everywhere with respect to an invariant measures, and with measurable  families of smooth manifolds, which coincide if they intersect and which fill  a set of full measure. 

We will denote by $\wl_\am^-(x)$ and $\w_\am^-(x)$ correspondingly the local and global
stable manifolds for the diffeomorphism $\am$ at a  point $x$ regular with respect to that
diffeomorphism. The global manifold is an immersed Euclidean space and is defined uniquely.  
Any local manifold is a $C^{1+\epsilon}$ embedded open disc in a Euclidean space. Its germ 
at $x$ is uniquely defined and  for any two choices  their intersection is an open neighborhood 
of the point $x$ in each of them. On a compact set of arbitrarily large measure, called a 
{\em Pesin set}, the local stable manifolds can be chosen of a uniform size and depending 
continuously in the $C^{1+\epsilon}$ topology. 

The local and global unstable manifolds $\wl_\am^+(x)$ and $\w_\am^+(x)$
are defined as the stable manifolds for the inverse map $\a(-\m)$.

Recall that the stable and unstable manifolds $\w_\am^-$ and $\w_\am^+$ are tangent 
to the  (almost everywhere defined)  stable and unstable distributions $E_\am^-$ and 
$E_\am^+$ accordingly. 
These distributions are the sums of the  distributions corresponding to the negative and 
positive Lyapunov exponents for $\am$ respectively. At the moment, we do not know 
the dimensions of those 
distributions.  We only know  from $(\mathcal E)$ that both distributions are 
non-trivial since by the Ruelle inequality  positive entropy implies existence of both 
positive and negative Lyapunov exponents.  

Corresponding stable and unstable  manifolds for the linear action $\ao$  will be 
denoted by the same symbols without the tilde $\, \tilde{}$. Of course those manifolds 
are affine, and they are defined everywhere, not just on large sets as for the nonlinear 
action $\a$. 

\subsection{Preservation of Weyl chambers under the semiconjugacy}\label{weyl}

For the following two lemmas we do not  need to assume that the  linear action $\ao$ 
is Cartan. 
It is sufficient  to assume that every element of $\ao$ other than identity is 
hyperbolic.\footnote{Even this condition may be weakened}

\begin{lemma}  \label{stable to stable}
For any element $\m\in \Zk\setminus\{ 0\}$ the following inclusions hold
  $$ 
       h(\w _\am ^- (x)) \subset \wo _\aom ^- (hx) \quad \text{ and } \quad
       h(\w _\am ^+ (x)) \subset \wo _\aom ^+ (hx),
  $$
   $$ 
       h(\wl _\am ^- (x)) \subset \wol _\aom ^- (hx) \quad \text{ and } \quad
       h(\wl _\am ^+ (x)) \subset \wol _\aom ^+ (hx).
  $$
on the set of full measure $\mu$ where $\w _\am ^- (x)$ and $\w _\am ^+ (x)$ 
exist.
\end{lemma}

\proof
The global stable manifold $\w _\am ^- (x)$ is the set of all points $y \in \tk$
for which $\dist (\a (n \m) x, \a (n \m)y) \to 0$ as $n \to \infty$. By continuity
of $h$ this implies that $\dist (\ao (n \m) hx, \ao (n \m)hy) \to 0$ as $n \to \infty$,
and hence $hy$ belongs to the stable manifold $ \wo _\aom ^- (hx)$. 

If $y$ is a point in the local stable manifold $\wl _\am ^- (x)$, it follows that the 
distance $\dist (\a (n \m) x, \a (n \m)y)$ remains small for all $n >0$. By continuity
of $h$, the distance  $\dist (\ao (n \m) hx, \ao (n \m)h)$ also remains small for all $n >0$. 
Since $\aom$ is uniformly hyperbolic, this implies that $hy$ belongs to the local stable 
manifold $\wol _\aom ^- (hx)$. 

The corresponding statements for the unstable manifolds follow by taking inverses.
$\QED$
\vskip0.3cm

\begin{lemma}  \label{Weyl chambers}
The Lyapunov half-spaces and Weyl chambers for $\a$ with respect to the measure  
$\mu$ are the same as the Lyapunov half-spaces and Weyl chambers for $\ao$.
Hence the Lyapunov exponents for $\a$ can be numbered $\tilde\chi_i,\,\, i=1,\dots,k+1$
in such a way that  $\tilde\chi_i=c_i\chi_i$ where $c_i$ is a positive scalar. 
\end{lemma}

\proof
Suppose that a Lyapunov hyperplane $L$ of $\ao$ is not a Lyapunov hyperplane 
of $\a$. Then there exist $\m,\n \in \Zk$ which lie  on the  opposite 
sides of $L$ so that $\wl ^-_\am = \wl ^-_\an$ but $\wol ^-_\aom \not= \wol ^-_\aon$.

Let $\La$ be the intersection of a Pesin set for $\am$ with a Pesin set for $\an$.
Consider  a  point $x\in \La$ such that any open neighborhood of $x$  intersects $\La$
by a set of positive measure $\mu$. By the previous lemma we have
 $$
   h(\wl ^-_\am(x))=h(\wl ^-_\an(x)) \subset (\wol ^-_\aom(hx) \cap  \wol ^-_\aon(hx)) 
   $$
   and
   $$
   h(\wl ^+_\am(x)) \subset \wol ^+_\aom(hx).
 $$
Let $R$ be the intersection of $\La$ with a neighborhood of $x$ 
sufficiently small compared to the size of the local manifolds at points of $\La$.
Then for any point $y \in R$ the intersection $\wl ^+_\am (x) \cap \wl^-_\am (y)$
consists of a single point $z_1$. Similarly, $\wl ^-_\am (x) \cap \wl^+_\am (y) =\{ z_2\}$
and hence $\wl ^-_\am (z_1) \cap \wl^+_\am (z_2) =\{ y\}$. By the previous lemma the
latter implies that $\wol ^-_\aom (hz_1) \cap \wol^+_\aom (hz_2) =\{ hy\}$. Using the
inclusions above we see that the image $h(R)$ is contained in the direct product 
$V =(\wol ^-_\aom (hx) \cap \wol ^-_\aon (hx) ) \times \wol ^+_\aom (hx) $. 
Since $\wol ^-_\aom \not= \wol ^-_\aon$, we conclude that  $V$
is contained in a subspace of dimension at most $k$. Hence $\la (V)=0$ which 
contradicts the fact that $\la (h(R)) \ge \mu(R) >0$. We conclude that any Lyapunov 
hyperplane of $\ao$ is also a Lyapunov hyperplane of $\a$. Recall that $\ao$ is
Cartan and thus have the maximal possible number, $k+1$, of Lyapunov hyperplanes.
Hence $\a$ also has exactly $k+1$ distinct Lyapunov hyperplanes, which coincide with
the Lyapunov hyperplanes of $\ao$. In particular, all Lyapunov exponents of $\a$ do
not vanish.

It follows that for either action there is exactly one Lyapunov exponent that corresponds
to a given Lyapunov hyperplane. It remains to check that for every Lyapunov hyperplane 
$L$ the  corresponding Lyapunov exponents of $\a$ and $\ao$ are {\em positively} 
proportional. Suppose that for some $L$ the corresponding Lyapunov exponents are 
negatively proportional. Let $W$ be the corresponding Lyapunov foliation for $\ao$.
We can take $\m$ close to $L$ in the negative half-space 
of the corresponding Lyapunov exponent for $\a$ and $\n$ sufficiently close to $\m$ across $L$
in the negative half-space for $\ao$, so that $\m$ and $\n$ are not separated from $L$
by any other Lyapunov hyperplane. Then we observe that $\wl _\am ^+ \subset \wl _\an ^+$ 
and that $W$ is contained neither in $\wol _\aom ^-$ nor in $\wol _\aon ^+$.

We choose $\La$, $x$, and $R$ as above. Using Lemma \ref{stable to stable} 
we obtain
  $$
     h(\wl _\am ^- (x)) \subset \wol _\aom ^- (hx) \quad \text{and} 
  $$
  $$
    h(\wl _\am ^+ (x)) \subset h(\wl _\an ^+ (x)) \subset \wol _\aon ^+ (hx).
  $$
As above, these inclusions imply that the image $h(R)$ is contained in 
$V=\wol _\aom ^- (hx) \times \wol _\aon ^+ (hx)$. We observe that $V$ lies in a 
subspace that does not contain $W(hx)$ and thus has dimension at most $k$. 
This again contradicts the fact that $\la (h(R)) \ge \mu(R) >0$.
$\QED$

Let us summarize the conclusions for the case of Cartan actions.

\begin{corollary}If $\ao$ is Cartan all  Lyapunov characteristic exponents for the action 
$\a$  with respect to measure $\mu$ are simple, no two of them are  proportional and  
the counterpart of property $(\mathcal C)$ holds.

For every Lyapunov exponent $\tilde\chi_i$ its Lyapunov distribution integrates to an 
invariant family of one-dimensional manifolds defined $\mu$ almost everywhere. This
family will be referred to as the Lyapunov foliation corresponding to $\tilde\chi_i$.
The semiconjugacy $h$ maps these local (corr. global)  manifolds to the local 
(corr. global) affine integral manifolds for the exponents $\chi_i$. 
\end{corollary}

%%%%%%%%%%%%%%%%  Proof of Main Theorem %%%%%%%%%%%%%%
%%%%%%%%%%%%%%%%  Proof of Main Theorem %%%%%%%%%%%%%%

\section{Proof of Theorem~\ref{thm-projection}}

Throughout this section we fix one of the Lyapunov exponents of $\a$. 
We denote 
by $L$ the corresponding Lyapunov hyperplane in $\Rk$,  
by $E$ the corresponding one-dimensional Lyapunov distribution, and 
by $\w$ the corresponding Lyapunov foliation.  
Then $\w$ is the one-dimensional stable foliation for some element $\am$, $\m \in \Zk$.
The notions of regularity and Pesin sets will refer to the corresponding notions 
for such an element.

In this section we study  properties of the action $\a$ related to $\w$. 
We will show that the conditional measure $\mw _x$ on the leaf $\w (x)$ is 
absolutely continuous for $\mu$ almost every $x$. We then conclude the proof
of Theorem~\ref{thm-projection} by showing that the absolute continuity of $\mu$
follows from the absolute continuity of conditional measures for every Lyapunov 
foliation.

%%%%%%%%%%%%%%%%  Affine Structures %%%%%%%%%%%%%%

\subsection{Invariant affine structures on leaves of Lyapunov foliations}
The following proposition gives a family of $\a$-invariant affine parameters 
on the leaves of the Lyapunov foliation $\w$. By an affine parameter we 
mean an atlas with affine transition maps. 

\begin{proposition}       \label{affine structures} 
There exists a unique 
measurable family of $\Ce$ smooth $\a$-invariant affine 
parameters on the leaves $\w (x)$. Moreover, they depend uniformly 
continuously in $\Ce$ topology on $x$ within a given Pesin set.
\end{proposition}

\begin{remark} Notice that those transition maps may  not always  preserve orientation.
In fact in some situations measurable choice of orientation is not possible. This however is completely irrelevant  for our uses of affine structures.
\end{remark}

  \begin{remark} In the proof below we use  the counterpart of the property 
 $(\mathcal C)$ for $\a$  and do not use  existence of a semiconjugacy with $\ao$. In fact, the assertion is true under a more general condition. 
 Namely, let $\chi$ be a simple (multiplicity one) Lyapunov exponent for an ergodic 
hyperbolic measure $\mu$ for a $\Ce$ diffeomorphism with an extra condition 
that there are no other exponents proportional to $\chi$ with the coefficient of 
proportionality greater than one. Then the  Lyapunov distribution for $\chi$ is 
integrable $\mu$ almost everywhere  to an invariant family of one-dimensional 
manifolds  and invariant affine parameters still exist. 

In the $C^2$ case one-dimensionality of the Lyapunov foliation may be  replaced by 
the following {\em bunching condition}: Lyapunov exponents may be positively proportional  with coefficients of proportionality  between 1/2 and 2. The coarse Lyapunov 
distribution  is always integrable and in this case the integral manifolds  admit  a unique invariant family of smooth affine structures. 

Proofs of these statements  can be obtained using non-uniform versions of the methods
from \cite{G}.
 \end{remark}

\proof The proposition is established using three lemmas below. We take
an element $\m \in \Zk$ such that $\w$ is the stable foliation of $\am$.
Then we apply Lemma \ref{linearization} with $f=\am$ to obtain the family
$H$ of non-stationary linearizations. Lemma \ref{consistency} then shows
that these non-stationary linearizations give an affine atlas. Since the 
linearization $H$ is unique by Lemma \ref{linearization unique}, the family 
$H$ linearizes any diffeomorphism  which commutes with $f$. Indeed, 
if $g\circ f=f\circ g$, then it is easy to see that $dg^{-1} \circ H_{g(\cdot)} \circ g\,$ also 
gives a non-stationary linearization for $f$, and hence $\,H \circ g=dg \circ H$.
Therefore, $H$ provides a non-stationary linearization for every element of the 
action $\a$, i.e. the action is affine with respect to the parameter. 
$\QED$
\vskip0.3cm

\begin{lemma}       \label{linearization} 
Let $\w$ be the one-dimensional stable foliation of a $C^{1+\epsilon}$ 
non-uniformly hyperbolic diffeomorphism $f$. Then for $\mu$ almost every point 
$x\in M$ there exists a $\Ce$ diffeomorphism $H_x: \w(x) \to E(x)=T_x \w$ 
such that
$$
  \begin{aligned}
  &(i) \;\;\; H_{fx}\circ f=D f \circ H_x, \\ 
  &(ii)\;\; H_x(x)=0 \text{ and } D_x H_x \text{ is the identity map,} \\
  &(iii)\; H_x \text{ depends continuously on } x 
     \text{ in }\Ce \text{ topology on a Pesin set}.\hskip3cm
  \end{aligned}
 $$ 
 \end{lemma}
 \vskip0.3cm
 
\proof 
We denote by $E$ the one-dimensional stable distribution for $f$.
We fix some background Riemannian metric $g$ on $M$ and denote
$$
 \df (x)  = \| Df (v) \|_{f x}  \cdot \|v\|_x^{-1}
$$
where $v \in E(x)$ and $\|.\|_x$ is the norm given by $g$ at $x$.

We first construct the diffeomorphism $H_x$ on the local manifold $\wl (x)$ as follows. 
Since $E(x)$ is one-dimensional, $H_x(y)$ for $y \in \wl (x)$ 
can be specified by its distance to $0$ with respect to the Euclidian metric on $E(x)$ 
induced by $g$. We define this distance  by integrating a H\"older continuous density  
\begin{equation}\label{Hx}
 |H_x (y)| = \int _x ^y \rho_x (z) dz 
 \end{equation}

where 
$$\rho_x (z) =\lim_{n \to \infty} \frac{\dfn (z)}{\dfn (x)} = 
    \prod_{k=0}^\infty \frac{\df (f^k z)}{\df (f^k x)} $$

For any point $z$ in the local manifold $\wl (x)$ we have 
$\dist (f^k z, f^k x) \le C(x) e^{-k\la} \dist (z,x)$
for all $k>0$. In particular, $f^k z$ remains in the local manifold $\wl (f^k x)$ 
even though the size of $\wl (f^k x)$ may decrease with $k$ at a slow exponential 
rate. The tangent space $E(s)=T_s \wl (f^k x)$ depends H\"older continuously on 
$s \in \wl (f^k x)$, with H\"older exponent $\epsilon$ and a constant which may increase 
with $k$ at a slow exponential rate. Since $f$ is $\Ce$, the same holds for $\df (z)$. 
We conclude that
$$
\left| \frac{\df (f^k z)}{\df (f^k x)}-1 \right| \le C(x) \dist (z,x) e^{-k(\la+\delta)}.
$$
This implies that the infinite product which defines $\rho_x (z)$ converges, and
that $\rho_x $ is H\"older continuous on $\wl (x)$. Moreover, the convergence
is uniform when $x$ is in a given Pesin set. Hence $\rho_x $ depends  
continuously in $C^\e$ topology on $x$ within a given Pesin set. Since $\rho_x (x)=1$, 
we conclude that \eqref{Hx} defines a $\Ce$ diffeomorphism satisfying conditions (ii) 
and (iii). To verify condition (i) we differentiate $H_{fx} (f(y)) = D_x f (H_x(y))$ with 
respect to $y$ and obtain $\rho_{fx} (fy) \cdot \df(y) = \df (x) \cdot \rho_x (y) $.
Since the latter is satisfied by the definition of $\rho$, the condition (i) follows
by integration.

Since $f$ contracts $\w$, we can extend $H$ to the global stable manifolds $\w (x)$ 
as follows. For $y \in \w (x)$ there exists $n$ such that $f^n(y)\in \wl (f^nx)$ and we
can set  
 $$
  H_x(y)=Df^{-n}\circ H_{fx} \circ f^n(y).
 $$
This defines $H_x$ on an increasing sequence of balls exhausting $\w (x)$ with
conditions (i) and (ii) satisfied by the construction. Condition (iii) is satisfied in the
following sense. $H_x$ is a $\Ce$ diffeomorphism with locally H\"older derivative.
Its restriction to a ball of fixed radius in $\w (x)$ centered at $x$ depends continuously
in $\Ce$ topology on $x$ within a given Pesin set.   $\QED$

\begin{remark}
In general, the regularity of the density $\rho_x$ on $\w (x)$ is the same 
as the regularity of the differential $Df$, and hence the function $H_x$ 
is as regular as $f$.
\end{remark}

\begin{lemma} \label{consistency}
Under the assumptions of Lemma \ref{linearization}, the map
$$
    H_y \circ  H_x^{-1}:\; E(x) \to E(y)
$$
is affine for any $x$ and $y$ on the same leaf of $\w$. Hence the 
non-stationary linearization $H$ defines affine 
parameters on the leaves of $\w$. 
\end{lemma}

\proof
By invariance under $f$, it suffices to consider $x$ and $y$ close, and
show that the map is affine in a neighborhood of zero.
We will show that the differential $D\left( H_y \circ  H_x^{-1}\right)$
is constant on $E(x)$. Consider $z \in \w(x)$ close to $x$ and $y$ and 
let $\bar z = H_x (z)$. From the definition of $H$ we have
$D_z (H_x)=\rho_x (z)$ and $D_z (H_y)(z)=\rho_y (z)$. Hence,
using the definition of $\rho$, we obtain
 $$
 D_{\bar z}\left( H_y \circ  H_x^{-1}\right)= D_{z}(H_y) \cdot D_{\bar z}(H_x^{-1})
 = D_{z}(H_y) \cdot \left(D_{z}(H_x) \right) ^{-1} =
 $$
 $$
 =\frac{\rho_y (z)}{\rho_x (z)}= \prod_{k=0}^\infty \frac{\df (f^k z)}{\df (f^k y)}
 \cdot \left( \prod_{k=0}^\infty \frac{\df (f^k z)}{\df (f^k x)} \right) ^{-1} =
 \prod_{k=0}^\infty \frac{\df (f^k y)}{\df (f^k x)}.
 $$
We conclude that the differential $D_{\bar z}\left( H_y \circ  H_x^{-1}\right)$ 
is independent of $\bar z$ and thus the map $H_y \circ  H_x^{-1}$ is affine. 
$\QED$
\vskip.3cm

\begin{lemma}       \label{linearization unique} 
The family of diffeomorphisms $\{ H_x \}$ satisfying conditions (i)-(iii) of 
Lemma \ref{linearization} is unique.
\end{lemma}
\vskip0.3cm

\proof
We note that it is sufficient for the proof to have $H_x$ defined only locally,
in a neighborhood of $x$ in $\wl (x)$.

Suppose that $H_1$ and $H_2$ are two families of maps 
satisfying (i)-(iii). Then the family of maps $G=H_1\circ H_2^{-1} : E \to E$ 
satisfies $G_{fx} \circ D_x f=D_x f \circ G_x$, and hence
 $$
  G_x=(D_x f)^{-1}\circ G_{fx} \circ D_x f= \dots=
      (D_x f ^n)^{-1}\circ G_{f^n x} \circ D_x f^n.
 $$
or, since $E$ is one--dimensional,
$$
	G_x(t)=(\dfn (x))^{-1} G_{f^n x} (\dfn(x) \cdot t).
$$

Since $\dfn(x) \to 0$ and since $G_x$ depends continuously in 
$C^1$-topology on $x$ in a Pesin set, we obtain using returns to such
a set that 
$$\frac{G_x(\dfn(x) \cdot t)}{\dfn(x) \cdot t} \to G'_x(0)=1$$ 
and hence 
$$
	G_x(t)=\lim_{n\to \infty} t \cdot \frac{G_{f^n x}(\dfn(x) \cdot t)}
	{\dfn(x) \cdot t} =t.
$$
Thus $G_x$ is the identity, and $H_1=H_2$
$\QED$
\vskip0.3cm

%%%%%%%%%%%%%% END Affine Structures %%%%%%%%%%%%%%

\subsection{Uniform growth estimates along the walls of Weyl chambers}
In the rest of  this section we consider suspensions of the actions $\ao$ and $\a$. 
 According to our convention we will use the same  notations for the  
suspension actions and associated objects. 

We fix a Pesin set $\La$ and a small $r>0$.  For $x \in \La$ we denote by 
$\bwr (x)$ the ball (interval) in the inner metric of  $\w (x)$ of radius $r$ centered at $x$. 
%We assume that $r$ is small enough so that $\bwr (x)$ is contained in the local 
%manifold ??
An important corollary of the existence of affine parameters is the 
following estimate of derivatives along $\w$.

\begin{lemma}  \label{derivative estimate}
For a given Pesin set $\La$ and $r>0$ there exists a constant $C=C(\La,r)$ 
such that for any $x \in \La$ and $\t \in \Rk$ satisfying $\at x \in \La$ 
 $$
   C^{-1}  \|D (\at )|_{E(x)} \| \le \|D (\at )|_{E(y)} \|  \le C  \|D (\at )|_{E(x)} \| 
 $$
for any $y \in \bwr (x)$ satisfying $\at y \in \bwr (\at x)$.
\end{lemma}

\proof
We use the affine parameter on $\w (x)$ given by Proposition \ref{affine structures} 
with respect to which $\at$ has constant derivative. More precisely, using the linearization
$H$ along the leaves of $\w$ we can write 
 $$
    \at |_{\w (x)}=(H_{\at x})^{-1} \circ D\at|_{E(x)} \circ H_x,
  $$  
and hence
 $$
   D\at|_{E(y)}=( D_{\at y} H_{\at x})^{-1} \circ D\at|_{E(x)} \circ D_y H_x.
 $$
Since $H_z|_{\bwr (z)}$ depends continuously in $\Ce$ topology on $z$ in the 
Pesin set $\La$, both $\|D_t H_z \|$ and $\|(D_t H_{z})^{-1}\|$ are uniformly 
bounded above and away from $0$ for all $z \in \La$ and $t \in \bwr(z)$. 
Hence the norms of the first and last term in the right hand side are uniformly bounded 
and the lemma follows.  $\QED$
\medskip

We consider $x \in \La$ and the ball $\bwr (x) \subset \wl (x)$. The image $h(\bwr (x))$ is 
contained in $\wol (x)$. We denote by $m_r (x)$ the radius of the largest ball (interval) in 
$\wol (hx)$ that is centered at $hx$ and contained in $h(\bwr (x))$. Then $m_r$ 
is a measurable function on $\La$.

\begin{lemma} \label{h estimate}
For any Pesin set $\La$ and $r>0$ the function $m_r$ is positive almost
everywhere on $\La$. Hence for any $\e>0$ there exists $m>0$ and a set 
$\La_{r,m} \subset \La$ with $\mu (\La \setminus \La_{r,m}) < \e$ so that 
$m_r (x) \ge m$ for all $x \in \La_{r,m}$.
\end{lemma}

\proof
Let $x\in \La$ be such a point that intersection of $\La$ with any neighborhood 
of $x$ has positive measure. Let $\m \in \Zk$ be an element such that 
$\w =\w^-_\am$.  Let $R$ be the intersection of $\La$ with a  sufficiently small
neighborhood of $x$. If $m_r (x)=0$ then 
$h(\wl^-_\am (x)) = \{ hx \}$. This implies, as in Lemma \ref{Weyl chambers}, that 
the image $h(R)$ is contained in $\wol^+_\aom (hx)$. But this implies that $\la (h(R))=0$, 
which is impossible since $\la (h(R) )\ge \mu(R) >0$. 
$\QED$
\vskip0.3cm

Using the derivative estimate in Lemma~\ref{derivative estimate} and  the 
topological semiconjugacy $h$ we obtain in  the next lemma  the crucial 
estimate for the derivatives of the elements in the Lyapunov hyperplanes.

We fix a Pesin set $\La$, $r>0$, and a set $\La_{r,m}$ as in Lemma \ref{h estimate}.

\begin{lemma}  \label{neutral estimate}
For a given set $\La_{r,m}$ there exists a constant $K$ such that
for any $\t$ in the Lyapunov hyperplane $L$
 $$
    K^{-1} \le \|D (\at )|_{E(x)} \|  \le K 
 $$
if both $x \in \La_{r,m}$ and  $\at x \in \La_{r,m}$.
\end{lemma}

\proof
First we note that it suffices to establish the lower estimate, then the upper 
estimate follows by applying it to $\a (-\t)$.

By uniform continuity of the semiconjugacy $h$ there exists $\delta >0$ such that 
for any $x$ the image $h(\tilde B_\delta (x))$ is contained in the ball 
$B_{m/2}(hx)$ in $\wol (x)$. By the choice of $\La_{r,m}$
we also have $\bwo_m (hx) \subset h(\bwr (x))$. Since $\t\in L$, $\aot$ is an isometry 
on $\wo$, and hence $\aot (\bwo _m (hx) ) = \bwo _m (\aot(hx))$. Then since $h$ is a 
semiconjugacy we obtain 
\begin{equation}
   \bwo _m (\aot(hx)) \subset (\aot \circ h) \, (\bwr (x))  = (h \circ \at) \, (\bwr (x)).
\label{3.9}    
\end{equation}
Together with the uniform continuity of $h$ this implies that $\at  (\bwr (x))$ cannot
be contained in $\tilde B_\delta (\at x)$. Indeed, otherwise we would have
$\bwo _m (\aot(hx)) \subset h (\tilde B_\delta (\at x)) 
\subset \at  (\bwr (x)) \subset B_{m/2}(\aot (hx)).$

Hence there exists $z \in \bwr (x)$ with $\dist(\at x, \at z)=\delta$. We may assume 
that $\delta < r$ and $z$ is chosen so that $\dist(\at x, \at y) < \delta$ for all 
$y\in \bwr (x)$ between $x$ and $z$. Using Lemma \ref{derivative estimate} we obtain  
  $$
  \delta =\dist(\at x, \at z) \le  \dist(x, z) \cdot \sup \|D (\at )|_{E(y)} \| <
  $$
  $$
  < r \cdot C \|D (\at )|_{E(x)} \| 
  $$
This implies that $\|D (\at )|_{E(x)} \|  > \frac{\delta}{Cr}$.    $\QED$

\subsection{Ergodicity along the walls of Weyl chambers}\label{ergodicitysingular}
 We will call an element $\t\in\Rk$ a {\em generic singular element} if it belongs to exactly one  Lyapunov hyperplane. The following lemma presents   a variation of an argument from \cite{KS3}
for  the present   setting.

\begin{lemma}  \label{transitivity on leaves}
Let $L$ be one of the Lyapunov hyperplanes in $\Rk$. Let $E$ and $\w$ be 
the corresponding Lyapunov distribution and foliation of $\a$.  
Then for any generic singular element $\t \in \Rk$ the corresponding partition 
$\xi_{\at}$ into the ergodic components of $\mu$ with respect to $\at$ is coarser than the measurable 
hull $\xi(\w)$ of the foliation $\w$.
\end{lemma}

\proof 
Consider a generic singular element $\t$ in $L$. Then the only non-trivial 
Lyapunov exponent that vanishes on $\t$ is the one with kernel $L$ and 
the corresponding Lyapunov distribution is $E$. Take a regular element 
$\s$ close to $\t$ for which this Lyapunov exponent is positive and all other 
non-trivial  exponents have the same signs as for $\t$. Thus 
$E_\as^+=E_\at^+\oplus E$ and $E_\as^-=E_\at^-$. 
Birkhoff averages with respect to $\at$ of any continuous 
function are  constant on the leaves of $\w^-_\at$. Since such averages 
generate the algebra of $\at$--invariant functions, we conclude that the
partition $\xi_\at$ into the ergodic components of $\at$ is coarser than 
$ \xi(\w_\at^-)$, the measurable hull of the foliation $\w^-_\at$. 
On the other hand, the measurable hulls $\xi(\w_\as^-)$ and $\xi(\w_\as^+)$ 
of both $\w_\as^-$ and $\w_\as^+$ coincide with the Pinsker algebra $\pi(\as)$. 
Since $\xi(\w _\as ^+)$ is coarser than $\xi (\w)$, we conclude
  $$
     \xi_\at \leq \xi (\w _\at^-) = \xi (\w_\as^-) = \pi(\as) =  \xi(\w _\as^+) \leq  \xi (\w).
         $$  $\QED$

%%%%%%%%%

\subsection{Invariance and absolute continuity of conditional measures}
Let $\w$ be one of the Lyapunov foliations of $\a$ (recall that it is one--dimensional), 
and let $L \subset \Rk$ be the 
corresponding Lyapunov hyperplane. We fix a Pesin set $\La$, $r>0$, and a 
set $\La_{r,m}$ as in Lemma \ref{h estimate}.

\begin{lemma}  \label{transitive group}
For $\mu$- a.e. $x \in \La_{r,m}$ and for $\mw _x$- a.e.  $y \in \La_{r,m} \cap \bwr (x)$
there exists an affine map $g : \w (x) \to \w (x)$ with $g (x)=y$ which preserves the 
conditional measure $\mw _x$  up to a positive scalar multiple Furthermore, the absolute
value of the derivative of this affine map is bounded away from zero and infinity uniformly 
in $x$ and $y$. The bounds depend on $r$ and $m$.
\end{lemma}

\proof We fix a generic singular element $\t \in L  \subset \Rk$. 
By Lemma \ref{transitivity on leaves}
the partition $\xi_{\at}$ into the ergodic components of $\mu$ for $\at$ is coarser 
than the measurable hull $\xi(\w)$ of the foliation $\w$.
Then there is a set $X_1$ of full $\mu$-measure such that for any $x \in X_1$
the ergodic component $E_x$ of $\at$ passing through $x$ is well-defined and 
contains $\w (x)$ up to a set of  $\mw _x$-measure 0. Let $\mu _x$ be the measure 
induced by $\mu$ on  $E_x$. 

%We fix $\s \in \Rk$ for which $\as$ expands $\w$. 
For $n>0$ we denote by $B^n (x)$ the image under $H_x^{-1}$ of the ball in
$T_x \w$ of radius $n$ centered at $0$, where $H_x$ comes from Lemma~\ref{linearization}. We note that the sets $B^n (x)$
 exhaust $\w (x)$, i.e. $\w (x) = \bigcup_{n>0} B^n (x)$. For almost every $x$ we can  normalize $\mw _x$ so that $\mw _x (B^n (x))=1$ and denote its restriction 
to $B^n (x)$ by $\mu _x^n$.
 
We use a fixed Riemannian metric to identify $T_x \w$ with $\R$ and then use $H_x$
to identify $B^n (x)$ with the interval $[-n,n]$. Thus we can consider the system 
of normalized conditional measures $\mu _{x}^n$ as a measurable function from 
the suspension manifold $M$ to the weak* compact set of  Borel probability 
measures on the interval $[-n,n]$. By  Luzin's theorem, we can take an increasing 
sequence of closed sets $K_i$ contained in the support of $\mu$ such that

\begin{enumerate}
\item  $\mu (K) = 1$, where $K=\bigcup \limits_{i=1}^{\infty} K_i$
\item  $\mu _{x}^n$ depends continuously on $ x \in K_i$ with respect
          to the weak$^*$ topology.
\end{enumerate}

Set  $X_2 = X_1 \cap K$. Since  by definition the transformation  $\at$ restricted to the ergodic  component 
$E_x$ is ergodic, the transformation induced by $\at$ on 
$X_1 \cap E_x \cap K_i \cap \La_{r,m}$ is also ergodic for any $i$. Hence the 
set $X_3$, which consists of points $x \in X_2$ whose orbit 
$\{ \a (m\t) \,x\} _{m \in \Z}$ is dense in a subset of full $\mu_x$ measure of  $X_1 \cap E_x \cap K_i \cap \La_{r,m}$ for all $i$,
has full measure $\mu$.

Let $x \in X_3 \cap \La_{r,m}$ and $y \in X_3  \cap \La_{r,m} \cap \bwr(x)$. Then 
$x,y \in X_1 \cap E_x \cap K_i \cap \La_{r,m}$ for some $i$. Hence there exists 
a sequence $m_k \rightarrow \infty$ such that the points 
$y_k = \a (m _k \t) \,x \in X_1  \cap E_x \cap K_i \cap \La_{r,m}$ converge to $y$. 
Let us consider the map 
  $$ 
     \phi_k = \a (m _k \t)|_{\w(x)} : \w(x) \to \w (y_k).
  $$ 
Since $x$ and $y_k= \a (m _k \t) \,x$ are both in $\La_{r,m}$, Lemma 
\ref{neutral estimate} shows that $K^{-1} \le \|D_x \phi_k \| \le K$  for all $k$.
The map $\phi_k$ is affine with respect to the affine parameters on $\w (x)$
and $\w (y_k)$. By Proposition \ref{affine structures}, the affine parameters depend 
continuously in $\Ce$ topology on a point in the Pesin set $\La$. Thus the affine 
parameters at $y_k$ converge to the affine parameter at $y$ uniformly on compact
sets in the leaves. Hence, by taking a subsequence if necessary, we may assume 
that $\phi_k$ converge uniformly on compact sets to an affine map 
$g_n : \w(x) \to \w(x)$ with $g_n (x) = y$. 

Since both $(\phi_k)_* \mu ^n _x$ and $\mu^n _{y_k}$ are conditional measures 
on the same leaf $\w(y_k)$, there exists a constant $c(k)>0$ such that 
  $$
  \mu ^n _{y_k} (\phi_k A) = c(k)  \mu ^n _x (A) \quad \text{ for any } \quad
   A \subset B^n(x) \cap \phi_k^{-1} (B^n_{y_k}).
  $$ 
Similarly, there exists a constant $c>0$ such that
 $$
   \mu ^n _{y} (A) = c  \mu^n _x (A) \quad \text{ for any } \quad
   A \subset B^n(x) \cap (B^n_{y}).
  $$ 
Since $\mu _x^n$ depends continuously on $x \in K_i$ with respect to the 
weak* topology, measures $ \mu _{y_k}^n$ weak* converge to the measure 
$\mu _y^n$.  Assuming that the boundary of $A$ relative to the leaf has
zero conditional measure, we obtain that
  $$
  c(k)  \mu^n _x (A) = \mu^n _{y_k} (\phi_k A) \to 
                                  \mu^n _{y}     (g_n A) =  c  \mu^n _x (g_n A)
  $$ 
and hence
 $$
   \mu^n _x (g_n A) = \frac{\lim c(k)}{c}  \mu^n _x (A) \quad \text{ for any } \quad
   A \subset B^n(x) \cap g_n^{-1} (B^n_{y}).
  $$ 
We obtain that $g_n$ preserves the conditional measure $\mw _x$ up to a scalar
on the set $C^n(x)=B^n(x) \cap g_n^{-1} (B^n_{y})$. We note that $C^n(x)$ contains 
$B^{n/K}(x)$ and also $\bwr (x)$, provided that $n$ is large enough. Since 
$\mw _x (\bwr (x))>0$, taking $A=C^n(x)$ we see that $\lim c(k)$ must be positive.

We conclude that for any $n>0$ there exists a set $X_4$ of full $\mu$
measure such that for any $x \in X_4  \cap \La_{r,m}$ and 
$y \in X_4  \cap \bwr(x)  \cap \La_{r,m}$ there exists an affine map $g_n$ of $\w(x)$ 
such that $g_n(x)=y$ and $g_n$ preserves $\mw _x$ up to a positive scalar on $C^n(x)$.

Repeating this construction for every $n>0$ we can choose a set $X$ of full 
measure $\mu$ such that for any $x \in X \cap \La_{r,m}$, 
$y \in X \cap \bwr (x) \cap \La_{r,m}$, and any $n$ there exists an affine map
$g_n : \w(x) \to \w(x)$ satisfying $g_n (x) = y$ and preserving $\mw _x$ up to 
a positive scalar on $C^n(x)$. We note that $\w (x) = \bigcup_{n>0} C^n (x)$. Hence taking 
a converging subsequence we obtain that for any $x \in X  \cap \La_{r,m}$ and 
$y \in X \cap \bwr(x)  \cap \La_{r,m}$ there exists an affine map $g$ of $\w(x)$ 
with $g(x)=y$ which preserves $\mw _x$ up to a positive scalar.  
This completes the proof of the lemma since we may assume that the set $X$ of 
full measure is chosen so that for $x \in X \cap \La_{r,m}$ the set $X \cap \bwr(x)$ 
has full $\mw _x$-measure.  $\QED$

\begin{lemma}  \label{conditional abs cont}
The conditional measures $\mw_x$  are absolutely 
continuous for $\mu$ - a.e. $x$. 
\end{lemma}

\proof
Let $A_x$ be the group of affine transformations of $\w (x)$, and let $G_x$ 
be the subgroup of $A_x$ consisting of elements which preserve $\mw _x$ 
up to a positive scalar multiple. 

We first observe that $G_x$ is a closed subgroup. Indeed, if
$g_n \to g$ in $A_x$ then $g_n (Z) \to g (Z)$ in Hausdorff metric for any
bounded  $Z \subset \w (x)$ and hence $\mw _x (g_n (Z)) \to \mw _x (g(Z))$ 
if the relative boundary of $g(Z)$ has zero conditional measure. This
implies that $(g_n)_\ast \mw _x \to g_\ast \mw _x$. We also have 
$(g_n)_\ast \mw _x = c_n \mw _x$, where $c_n = \mw _x (Z) /
\mw _x (g_n (Z))$ for any $Z$. Since $g$ is an invertible affine map we can 
choose $Z$ such  that with $\mw _x (Z)>0$, $\mw _x (g (Z))>0$, and 
$\mw _x (\partial (g (Z)))=0$. It follows that $c_n \to c = \mw _x (Z) /
\mw _x (g (Z)) >0$ and $g_\ast \mw _x = c\mw _x$.

Since any element $\at$ preserves affine parameters on the leaves of $\w$,
it maps the group $A_x$ isomorphically onto $A_{\at x}$. Since $\at$ also
preserves the conditional measures on the leaves of $\w$, it maps the 
subgroup $G_x$ isomorphically onto $G_{\at x}$ on the set of full measure 
$\mu$ where the conditional measures and affine parameters on the leaves 
of $\w$ are well defined.  Since  isomorphism  classes of closed subgroups  
of the group of affine transformations on the line form a separable space, 
ergodicity of $\at$ implies that the  groups $G_x$  are isomorphic  
$\mu$-almost everywhere. 

By Lemma \ref{transitive group}, for a given Pesin set $\La$ and for $\mw _x$ almost any 
$y,z \in \La_{r,m} \cap \bwr (x)$ there exists an affine map $g : \w (x) \to \w (x)$ 
preserving $\mw _x$ up to a scalar multiple with $g (y) =z$. Thus $G_x$ has an orbit of positive $\mw _x$ measure.
We note that the measures $\mw _x$ are non-atomic for $\mu$ almost every $x$,
otherwise the entropy would be zero for any element whose full unstable foliation 
is $\w$.   
Then it follows that $G_x$ can not be a discrete subgroup of $A_x$ for $\mu$ 
almost every $x$. Hence either $G_x=A_x$ or the connected component of the 
identity in  $G_x$ is a one-parameter subgroup of the same type on the set of full 
$\mu$ measure. Thus either  $G_x$ contains the subgroup of translations or it is 
conjugate to the subgroup of dilations.
\vskip0.3cm

\noindent (i) First  consider the case when $G_x$ contains the subgroup of translations.
For any $x$ and $y \in \w(x)$ we define  $c_x(y)$ by the equality $g \mw = c_x(y) \mw$, 
where $\mw$ is the conditional measure on $\w (x)$ and $g$  is a translation such that 
$g(x)=y$. Note that $c_x(y)$ is well defined. Indeed, such $g$ is unique, and the definition 
does not depend on a particular choice of $\mw$ since the conditional measures are defined 
up to a scalar multiple. We need to show that  $c_x(y) =1$ for all $y \in \w (x)$.

We note that $c_x(y)$ can be calculated as 
  $$c_x(y)=\dfrac{g_\ast \mw (A)}{\mw (A)}=\dfrac{\mw (g^{-1} A)}{\mw (A)}
        =\dfrac{\mw (A)}{\mw (g A)}
  $$
for any set $A$ of positive conditional measure.  Since we can take the test 
set $A$ such that the boundary of $g A$ relative the leaf has zero conditional 
measure, we conclude that for a fixed $x$ the coefficient 
$c_x(y)$ depends continuously on $y$.

We see that either for $\mu$- a.e. $x\;\; c_x(y) =1$ for all $y \in \w (x)$, or there 
exists a set $X$ of positive measure such that $c_x(y)$ is not identically equal 
to $1$ for $x \in X$. In the latter case for some $\epsilon >0$ we can define a 
finite positive measurable function 
$$\varphi_{\epsilon}(x)=\inf
\{r : \exists \; y \in \w (x) \; s.t.\; d(x,y)<r \; and \;
|c_x(y)-1|>\epsilon \}  $$
on some subset $Y \subset X$
of positive $\mu$-measure. By measurability there exists $N$ and a set
$Z$ of positive measure on which $\varphi_\epsilon$ takes values
in the interval $(1/N,N)$.  We will show that 
\begin{equation}\varphi_{\epsilon}(\a (n \t) x) \rightarrow 0 \text{ as } n \rightarrow \infty
\label{(3.3.1)} 
\end{equation} 
uniformly on $Z$ for an element $\t$ such that $\at$ contracts the foliation $\w$. 
Since this contradicts to the recurrence of the set  $Z$ we conclude that $c_x(y)$ must be identically equal to $1$. 

We will prove now that 
  \begin{equation}
     c_x(y) = c_{\at x}(\at y) \label{(3.3.2)} 
\end{equation}
for $\mu$-a.e. $x$ and $y \in \w (x)$.  Since the iterates of $\at$ 
exponentially contract the leaves of $\w$, this invariance property 
implies that $\varphi_{\epsilon}(\a (n \t) x) \leq C \lambda ^n \varphi_{\epsilon}(x)$, 
for some $C$,$\lambda >0$, hence \eqref{(3.3.2)} implies \eqref{(3.3.1)}.

To prove \eqref{(3.3.2)} we consider translation 
 $$
    f=\at \circ g \circ \a (-\t) \in G_{\at x}
 $$
We observe that $f (\at x) = \at y$ since $gx =y$. Hence we obtain
 $$
     c_{\at x}(\at y) = \dfrac{\mw  (B)}{\mw (f B)}   
 $$
for any set $B \subset \w(\at x)$ of positive conditional measure. Since 
$\at ( g A) = f (\at A)$ and since $\at _\ast \mw _x$ is a conditional measure 
on the leaf $\w(\at x)$ we obtain using $B=\at A$  as the test set that
 $$
    c_x(y) = \dfrac{\mw_x (A)}{\mw_x (g A)}=
                  \dfrac{(\at_\ast \mw _x) (\at A)}{(\at_\ast \mw_x)(\at ( g A))}= c_{\at x}(\at y).
 $$
\vskip0.3cm

\noindent({ii)  Now suppose that $G_x$ is conjugate to the subgroup of  dilations. 
In this case $G_x$ has a fixed point $0_x$ and acts simply transitively 
on each connected  component of $\w (x) \setminus \{ 0_x \}$. For any $x$ and 
$y$ in the same component we consider
  $$
      c_x(y)    = \dfrac{Jg \cdot \mw (A)}{\mw (g A)}
  $$
where $g \in G_x$ is such that $g(x)=y$ and $Jg$ is the absolute value of
the Jacobian with respect to the affine parameter. To show that measure $\mw _x$ 
is Haar it is sufficient to prove that for any $g \in G_x$
 \begin{equation}\label{equivariance}
 g_\ast \mw _x = Jg \cdot \mw _x
 \end{equation} 
For that it suffices to show that $c_x(y)=1$ 
identically on $\w (x)$ for $\mu$ almost every $x$. This can be established by
repeating the argument of the previous case. The only difference is that to prove 
\eqref{(3.3.2)} we need to note that for the map
$$
    f=\at \circ g \circ \a (-\t) \in G_{\at x}
 $$
we have $Jf=Jg$.  $\QED$

Notice  that at the end  we proved that $G_x=A_x$ for almost every $x$. 
\vskip 0.3 cm

\subsection{Conclusion of the proof}
In order to prove that $\mu$ is an absolutely continuous measure
it is sufficient to show that for a certain element  $\a(\m)$
\smallskip

($\mathcal P$) {\em Entropy $\h_{\mu}(\a(\m))$ is equal  both to the sum 
of the positive Lyapunov exponents and  to the absolute value of the sum 
of the negative Lyapunov exponents.} (See \cite{L, LY}).
\smallskip

First recall that  there are $2^{k+1}-2$  Weyl chambers for $\ao$ and any 
combination of positive and negative signs for the Lyapunov exponents, except 
for all  positive or all negative, appears in one of the Weyl chambers.  The same is 
true for $\a$ by Lemma~\ref{Weyl chambers}.
Denote the Lyapunov exponents  for $\a$ by $\chi_1,\dots, \chi_{k+1}$. 
Let $\mathcal C_i,\,\,i=1,\dots, k+1$, be the Weyl chamber on which the $\chi_i>0$ 
and $\chi_j<0$ for all $j\neq i$.
Notice that we use notations different from those of Section~\ref{section-prelim}.

Consider $\m\in\mathcal C_i$. Since the conditional measure on $\w ^+_{\am}$
is absolutely continuous by Lemma~\ref{conditional abs cont}, we obtain that
$$\h_{\mu}(\a(\m))=\chi_i(\m)$$
for any $\m\in\mathcal C_i$.  By the Ruelle entropy inequality
$\h_{\mu}(\a(\m)) \le -\sum_{j\neq i}\chi_j(\m)$ and hence
$$\sum_{j=1}^{k+1}\chi_j(\m)\le0.$$
If $\sum_{j=1}^{k+1}\chi_j(\m)=0$ then ($\mathcal P$) holds 
and the proof is finished. 

Thus we have to consider the case when $\sum_{j=1}^{k+1}\chi_j(\m)<0$ for all $\m$ 
in all Weyl chambers $\mathcal C_i,\,\,i=1,\dots k+1$. This implies that 
$\bigcup_{i=1}^{k+1}\mathcal C_i$ lies in a negative half space of the linear functional
$\sum_{j=1}^{k+1}\chi_j$. But this is impossible since there exist elements 
$\t_i\in \mathcal C_i,\,\,i=1,\dots k+1$  such that $\sum_{i=1}^{k+1}\t_i=0$.
$\QED$

\section{Proof of Theorems~\ref{thm-finitecover} and \ref{thm-exponents}}

\subsection{Rigidity of the expansion coefficients}
We consider the  suspension action of $\a$. 
Let $\chi$ be one of the Lyapunov exponents of $\a$. Let $E$ be the 
corresponding Lyapunov distribution and $L=\ker \chi \subset \Rk$ be the corresponding 
Lyapunov hyperplane.

\begin{lemma} The restriction of $\a$ to $L$ is ergodic.
\end{lemma}

\proof By Lemma \ref{transitivity on leaves},
the  partition $\xi_{L}$ into the ergodic components of $\mu$ is coarser than the 
measurable hull $\xi(\w)$ of the foliation $\w$, which coincides with the Pinsker
algebra of a regular element in $\Rk$. Since we have established that $\mu$ is
absolutely continuous, the Pinsker algebra on $\tk$ of a regular element in $\Zk$
is at most finite \cite{P}. 

 Then on the suspension manifold $M$ the Pinsker algebra 
 is given by the corresponding 
finite partitions of the fibers of the suspension. Since $L$ is an irrational hyperplane
in $\Rk$, its action on $\T ^k$ in the base of the suspension is uniquely ergodic,
and hence $\xi_{L}$ is at most finite. Since the action $\a$ is ergodic,
 $\xi_{L}$ is trivial since the stationary subgroup in $\Rk$
of any $L$-invariant set  has to have finite index and hence  must coincide with $\Rk$.
 $\QED$

\begin{lemma}  \label{normalization}
There is a measurable metric on $E$ with respect to which 
\begin{equation} \label{constant scaling}
\|D \at v \| =e^{\chi(\t)} \|v\|
\end{equation}
for any $\t \in \Rk$, $\mu$-a.e. $x$, and any $v\in E(x)$. 
Such a measurable metric is unique up to a scalar multiple.  

\end{lemma}

\proof
First we construct a measurable metric $g$ on $E$ which is preserved by an 
ergodic element $\t \in L$. In other words, \eqref{constant scaling} is satisfied
with respect to $g$ for this element $\t$. Then we will show that such a metric 
is unique up to a scalar multiple. The uniqueness easily implies that 
\eqref{constant scaling} is satisfied for all $\t \in \Rk$.

Let $\La '=\La_{r,m}$ and constant $K$ be as in Lemma \ref{neutral estimate}. 
Ergodicity of $\a|_L$ implies that there exists an ergodic element $\t \in L$. 
We fix such an element $\t$, and let $X$ be an invariant set of full measure 
consisting of points whose $\at$ orbits visit $\La '$ with frequency 
$\mu (\La ')$. 

We fix some background Riemannian metric $g_0$ on $M$.  
We use the notations $\de _x \at = D (\at )|_{E(x)}$ and
$$
\| \de _x \at \| = \| \de_x \at (v) \|_{\at x}  \cdot \|v\|_x^{-1}
$$
where $v \in E(x)$ and $\|.\|_x$ is the norm given by $g_0$ at $x$.

We define a measurable renormalization function $\phi$ as follows. 
\begin{equation} \label{renorm}
\phi (x) = \sup \{ \|\de _x \a (n\t) \| : n \in \N,\;   \a (n\t)x \in \La ' \} 
\end{equation}
We note that by Lemma \ref{neutral estimate} the supremum is bounded by $K$ 
for any $x \in \La '$. More generally, the supremum is finite for any point whose 
$\at$ orbit visits $\La '$. Thus the function is well defined and finite on $X$. Using 
\eqref{renorm} we obtain
$$
\frac{\phi (\at x)}{\phi (x)} = 
\frac{\sup \{ \|\de _{\at x} \a (n\t) \| :  n \in \N,\;  \a ((n+1)\t)x \in \La ' \}}
       {\sup \{ \|\de _x \a (n\t) \| : n \in \N,\;  \a (n\t)x \in \La ' \}}
 $$
 $$      
=\frac{\sup \{ \|\de _{\at x} \a (n\t) \| :  n \in \N,\;  \a ((n+1)\t)x \in \La ' \}}
{\sup \{ \|\de _x \at \| \cdot \|\de _x \a (n\t) \| : n \in \N,\;  \a ((n+1)\t)x \in \La ' \}}      
$$

$$
=\|\de _x \at \|^{-1}
$$
This means that with respect to the renormalized Riemannian metric $g=\phi g_0$
we have 
$$
\|\de _x \at \|_g=\|\de _x \at \| \cdot \frac{\phi(\at x)}{\phi( x)}=1.
$$

Suppose that \eqref{constant scaling} is satisfied for the fixed $\t$ with respect
to another Riemannian metric $\psi g_0$ on $E$. Then equation \eqref{constant scaling}
implies that
$$
 \|\de _x \at \| \cdot \frac{\psi(\at x)}{\psi( x)}=\|\de _x \at \|_{\psi g_0} =1
 $$
 $$
 =\|\de _x \at \|_{\phi g_0}= \|\de _x \at \| \cdot \frac{\phi(\at x)}{\phi( x)}
$$
and hence
$$
\frac{\psi(\at x)}{\phi(\at x)}=\frac{\psi( x)}{\phi( x)}
$$
By ergodicity of $\at$ we conclude that $\psi=\kappa \phi$, where $\kappa$ 
is a constant.

For any other element $\s \in  \Rk$ consider the metric $\as_\ast g$. 
By commutativity, this metric is again preserved by $\at$. From the uniqueness
we obtain $\as _\ast g = \kappa(\s) \cdot g$ where $\kappa(\s)$ is a positive constant.  Let us show that 
$$\log\kappa(\s)=\chi (\s).$$ 
Indeed, let $\La$ be a set of positive
measure on which $C^{-1} < \phi  < C$ for some constant $C$. Since
$$
\kappa(\s) = \|\de _x \a (n\s) \|_g=\|\de _x \as \| \cdot \frac{\phi(\a (n\s) x)}{\phi( x)}
$$
we obtain 
\begin{equation}\label{est}
C^{-2}  \kappa^n(\s) < \|\de _x \a (n\s) \| < C^2 \kappa^n(\s)
\end{equation}
if both $x$ and $\a (n\s)x$ are in $\La$. By recurrence, for almost every $x\in\Lambda$ 
there exists a sequence  of natural numbers $n_i\to\infty$ such that $\a(n_i\s)\in \La$. 
Since for almost every $x$
$$\chi(\s) = \lim_{i\to\infty} n_i^{-1} \log\|\de _x \a(n_k\s)\|$$
we conclude using \eqref{est} that $\chi (\s)=\log \kappa(\s)$.  $\QED$
\vskip0.3cm

\subsection{Smoothness of semiconjugacy along Lyapunov foliations}

\begin{lemma}  \label{transitive group 2}
For almost every $x$ the semiconjugacy $h$ intertwines the 
actions of the groups of translations of $\w (x)$ and $\wo (hx)$. More precisely, for any
translation $\tilde \tau$ with respect to the affine structure on $\w (x)$ there
is a translation $\tau$ of $\wo (hx)$ with $h\circ \tilde \tau =\tau \circ h$.
\end{lemma}

\proof 
The proof of this lemma closely follows the proof of Lemma \ref{transitive group}.  
Let $\La$ be a Pesin set. By Lemma \ref{affine structures}, the affine parameters 
depend continuously in $\Ce$ on a point in $\La$.
  
By Luzin's theorem, the measurable metric from Lemma \ref{normalization} 
is uniformly continuous on sets of large measure. Hence we can take an increasing 
sequence of closed sets $K_i$ such that

\begin{enumerate}
\item  $\mu (K) = 1$, where $K=\bigcup \limits_{i=1}^{\infty} K_i$
\item  the measurable metric depends continuously on $ x \in K_i$.
\end{enumerate}

As in the previous lemma we fix an ergodic element $\t \in L$. Then the transformation 
induced by $\at$ on $K_i \cap \La$ is also ergodic for any $i$. 
Hence, there is an invariant full measure $\mu$ set $X\subset K$ of points $x$ 
whose orbit $\{ \a (m\t) \,x\} _{m \in \Z}$ is dense in $K_i \cap \La$ for all $i$.

Let $x \in X$ and $y \in \w (x) \cap X \cap \La$. Then $y \in  K_i \cap \La$ for some 
$i$. Hence there exists a sequence $m_k \rightarrow \infty$ such that the points 
$y_k = \a (m _k \t) \,x \in K_i \cap \La$ converge to $y$. 
Let us consider the  affine map 
  $$ 
     \phi_k = \a (m _k \t)|_{\w(x)} : \w(x) \to \w (y_k).
  $$ 
We normalize the affine parameters using the measurable metric. 
%constructed in Lemma \ref{normalization} 
Then  $\phi_k$ 
is an isometry with respect to the normalized parameters at $x$ and $y_k$. The
normalized parameters vary continuously on $K_i \cap \La$. Since $y$ and $y_k$ 
are both in $K_i \cap \La$, the normalized affine parameters at $y_k$ converge to 
the normalized affine parameter at $y$ uniformly on compact sets. Hence, by taking 
a subsequence if necessary, we may assume that $\phi_k$ converge to an isometry 
$g : \w(x) \to \w(x)$ with $g (x) = y$. We also note that $y_k \to y$ implies that 
$hy_k \to hy$, and the maps 
$$ 
     \psi_k = \ao (m _k \t)|_{ \wo(hx)} : \wo(hx) \to \wo (hy_k).
 $$ 
are isometries. By taking a subsequence if necessary, we may assume that 
$\psi_k$ converge to an isometry $f : \wo(hx) \to \wo(hx)$ with $f (hx) = hy$.
Since $h$ is a semiconjugacy we obtain $h\circ g  =f \circ h$.

Let $G_x$ be the set of all isometries $g$ of $\w (x)$ for which there exists an 
isometry $f$ of $\wo(x)$ with $h\circ g  =f \circ h$. It is easy to see that $G_x$
is a closed subgroup of the group of affine transformations of $\w (x)$.

Since a set of full measure can be exhausted by Pesin sets we obtain that 
for almost every point $x$ and for $\mw _x$ almost every $y \in \w (x)$
there exists an isometry $g_{xy}\in G_x$ with $g (x) = y$. We note that by Lemma 
\ref{conditional abs cont}, for almost every point $x$ the conditional measure 
$\mu_x$ is Haar with respect to the affine parameter. Hence we conclude that for 
almost every point $x$ there is a dense set of points $y \in \w (x)$ for which 
there exists an isometry $g_{xy} \in G_x$. Since $G_x$ is closed this implies that
$G_x$ acts transitively on $\w (x)$ and thus contains the subgroup ${\cal T} _x$ of 
translations of $\w (x)$. The corresponding isometries of $\w (x)$ also have to
be translations and the lemma follows.
 $\QED$

\begin{lemma}  \label{smoothness}
For almost every point $x$ and every Lyapunov foliation $\w$ the 
semiconjugacy $h$ is a $\Ce$ diffeomorphism from $\w (x)$ into $\wo (hx)$.
\end{lemma}

\proof 
This follows immediately from Lemma \ref{transitive group 2}. Indeed, the correspondence
$\tilde \tau \to \tau$ is a continuous isomorphism between the groups of translations 
$\tilde {\cal T}$ and ${\cal T}$ of $\w(x)$ and $\wo(x)$ respectively. Hence there exists
$a \in \R$ such that if $\tilde \tau (y) = y + t$ for $y \in \w(x)$ then $\tau (z) 
= z + at$ for $z \in \wo(hx)$. Then $h\circ \tilde \tau = \tau \circ h$ implies
that $h|_{\w{x}}$ is a linear map with respect to the affine parameter on $\w(x)$ and
the standard affine parameter on $\wo(x)$. Since the affine parameter on $\w(x)$
is given by a $\Ce$ diffeomorphism, then so is $h$.   $\QED$
 \vskip 0.3cm

 \subsection{Conclusion of the proof of Theorem~\ref{thm-exponents}}\label{exponents-proof}
 
Fix a Lyapunov  foliation $\w$. Let  $\Lambda$  be a set of positive measure  
such that  the semiconjugacy $h$ is differentiable at every $x\in\Lambda$ and 
the  derivative $L(x)$ of the $h$ along $\w$ and its inverse are both bounded  
by a constant $C$. Such a set exists by Lemma~\ref{smoothness}.  Suppose
that both $x$ and $\a(\t)(x)$ are in $\Lambda$. Let $v$ be a tangent vector at 
$x$ to $\w$. We have 
 \begin{equation}\label{almostconstant} 
   \|D(\a(\t))v\|=L(x)\|D(\ao(\t))|_{\wo(hx)}\| L^{-1}(\a(\t)(x)\|v\|.
 \end{equation}
Since $\ao$ is a linear action $\|D(\ao(\t))|\wo(hx)\|=\exp\chi(\t)$, where $\chi$ 
is the Lyapunov exponent of  $\ao$ corresponding to the  foliation $\wo$.
Since by assumption 
 $$
  C^{-1}<\min \{L(x),  L^{-1}(\a(\t)(x))\}<\max \{L(x),  L^{-1}(\a(\t)(x))\}< C
  $$ 
we obtain from \eqref{almostconstant}
\begin{equation}\label{doublebound}
  C^{-2}\exp\chi(\t)\|v\|<\|D(\a(\t))v\|<C^2\exp\chi(\t)\|v\|.
 \end{equation}
Now let $\tilde\chi$ be the Lyapunov exponent of the action $\a$ corresponding 
to the foliation $\w$. Take $\s\in\Rk$ such that $\ao(\s)$ is ergodic (the set of such 
$\s$ is dense in $\Rk$).  Then for almost every $x\in\Lambda$ one can find a 
sequence  of natural numbers $n_k\to\infty$ such that $\a(n_k\s)\in \La$. Since 
for almost every $x$ and for $v\in T_x\w$
$$\lim_{k\to\infty}\frac{\log\|D\a(n_k\s)(v)\|}{n_k}=\tilde\chi(\s)$$
we conclude from \eqref{doublebound} that $\tilde\chi(\s)=\chi(\s)$. 
Since this is true for a dense set of $\s$, this implies that $\tilde\chi=\chi$
$\QED$
\vskip 0.3 cm

\subsection{Conclusion of the proof of Theorem~\ref{thm-finitecover}}
For every Weyl chamber $\mathcal C_i$ we choose an element 
$\m _i \in \Zk \cap \mathcal C_i$. For every $\m _i$ we choose a Pesin set
$\La _i$ and let $\La = \bigcap_i \, \La _i$. For a point $x$ in $\La$ we denote
by $B_r (x)$ the ball in $\tk$ of radius $r$ centered at $x$. We fix $r$ sufficiently 
small compared to the size of the local manifolds at points of $\La$.

\begin{lemma}  \label{local injectivity}
The semiconjugacy $h$ is injective on $B_r (x) \cap \La$ for any $x\in\La$.
\end{lemma}

\proof 
Let $y$ be a point in $B_r (x) \cap \La$ different from $x$. Then there exists
an element $\m _i$ such that $\wl ^+_{\a (\m _i)}(x)$ is $k$-dimensional
and does not contain $y$. Indeed, the intersection of all $k$-dimensional
local unstable manifolds through $x$ contains only $x$ itself. We will denote  
in this proof $\wl ^+_{\a (\m _i)}$ by $\tilde F$ and the complimentary one-dimensional
local Lyapunov foliation $\wl ^-_{\a (\m _i)}$ by $\tilde W$. By Lemma \ref{stable to stable},
$\tilde F(z) \subset F(hz)$ and $\tilde W(z) \subset W(hz)$ for any $z \in B_r (x) \cap \La$,
where $F$ and $W$ are the corresponding local foliations for $\ao$. Since both $x$ 
and $y$ are in $B_r (x) \cap \La$, the intersection $\tilde W(x) \cap \tilde F(y)$ consists
of exactly one point $z \in B_r (x)$. Then $hz$ is the unique point in the intersection 
$W(hx) \cap F(hy)$. Suppose now that $hx=hy$. Then $hz = W(hx) \cap F(hx)$,
which means that $hz=hx$ and thus $h$ is not injective on $W(x)$. The latter, 
however, contradicts the fact that, according to Lemma \ref{smoothness}, $h$ is a 
diffeomorphism on $W(x)$.  $\QED$
 \vskip 0.3cm

Now we can complete the proof of Theorem~\ref{thm-finitecover} as follows.
By compactness, the set $\La$ can be covered by finitely many balls $B_r (x)$. 
Then the previous lemma implies that $h^{-1}(hx) \cap \La$ 
is finite for any $x \in \La$. Since we can choose the set $\La$ to have arbitrarily
large measure $\mu$, by taking an increasing sequence of such sets we can 
obtain an invariant set $A$ with $\mu (A)=1$ such that $h^{-1}(hx) \cap A$ is at 
most countable for all $x\in A$. Now we consider the measurable partition of $A$
into the preimages of points under $h$. Since the elements of this partition are
at most countable, the conditional measures are discrete. We see that $h|_A$
gives a countable extension of $(\ao, \lambda)$. Since $(\a, \mu)$ is ergodic
it follows that the the conditional measures of all points in the fiber are the same
and hence the extension is finite. $\QED$

%%%%%%%%%%%%%%%%% End Text %%%%%%%%%%%%%%%%%%

%%%%%%%%%%%%%%%%%%%%%%%%%%%%%%%%%%%%%%
%%%%%%%%%%%%%%%       Bibliography       %%%%%%%%%%%%
%%%%%%%%%%%%%%%%%%%%%%%%%%%%%%%%%%%%%%

{\small

\end{document}